\documentclass{amsart}
\usepackage{amssymb}
\usepackage{amsmath}
\usepackage{amsthm}
\usepackage[all]{xy}
\xyoption{2cell}

\newcommand{\Ol}{{\mathcal O}}

\newcommand{\proj}{\mathbb P}
\newcommand{\R}{\mathbb R}
\newcommand{\quadr}{\mathbb Q}
\newcommand{\pt}{{\mathbb P^3}}
\newcommand{\pd}{{\mathbb P^2}}
\newcommand{\qt}{{\mathbb Q^3}}
\newcommand{\pn}{{\mathbb P^n}}
\newcommand{\qn}{{\mathbb Q^n}}

\newcommand{\Z}{{\mathbb Z}}

\DeclareMathOperator{\cone}{NE}
\DeclareMathOperator{\pic}{Pic}

\newcommand{\rc}[2]{#1 \xymatrix{\ar@{-->}[r] & }{#2}}

\newtheorem{theorem}{Theorem}[section]
\newtheorem{lemma}[theorem]{Lemma}
\newtheorem{proposition}[theorem]{Proposition}

\newtheorem{teo}{Theorem}[subsection]

\newtheorem{prop}[teo]{Proposition}

\theoremstyle{definition}
\newtheorem{definition}[theorem]{\rm Definition}
\newtheorem{statement}[theorem]{\rm}

\theoremstyle{remark}
\newtheorem{remark}[theorem]{\rm Remark}

\theoremstyle{definition}

\theoremstyle{remark}

\numberwithin{equation}{section}

\begin{document}
\author{Carla Novelli}
\address{Dipartimento di Matematica Pura ed Applicata  \newline
\indent 
Universit\`a degli Studi di Padova \newline
\indent via Trieste, 63  \newline
\indent 
I-35121 Padova, Italy} 
\email{novelli@math.unipd.it}


\title{Extremal rays of non-integral $L$-length}

\begin{abstract}
Let $X$ be a smooth complex projective variety and let $L$ be a line bundle on it.
We describe the structure of the pre-polarized manifold $(X,L)$ for non integral values of the invariant $\tau_L(R):=-K_X\cdot\Gamma/(L \cdot \Gamma)$, where $\Gamma$ is a minimal curve of an extremal ray $R:=\mathbb R_+[\Gamma]$ on $X$ such that $L \cdot R>0$.
\end{abstract}

\maketitle

\section*{Introduction}
Let $X$ be a smooth complex projective variety and let $L$ be a line bundle on~$X$.
Assume that the canonical bundle $K_ X$ of $X$ is not nef.
In classical {\em Adjunction Theory}, $L$ is assumed to be ample;
so by the {\em Kawamata's Rationality Theorem} the invariant $\tau := \tau (X,L) = \mbox{min} \{t \in \R: K_X+t L \mbox{ is nef}\}$ is a positive rational number, called the {\em nefvalue} of $(X,L)$. 
By the {\em Kawamata--Shokurov Base Point Free Theorem} it is possible to consider the morphism defined by a (sufficiently large) multiple of the divisor $K_ X+\tau L$. 
Then the classification of polarized manifolds $(X,L)$ in~terms of the values of $\tau$ and based on the study of the structure of this morphism is a natural question.
The book \cite{bookBS} is a good reference for Adjunction~Theory.

One main obstruction to extending this study to the case when $L$ is merely nef is given by the possible existence of  cycles $Z \in \overline{\cone}(X)$ such that $K_X  \cdot Z<0$ and $L\cdot Z=0$. 
Clearly in this case the invariant $\tau$ is not defined. 
In \cite{BKLN-LPosRays} we circumvent this problem in the following way:
since $K_X$ is not nef, it is well-known that there exists (at least) an extremal ray $R:=\mathbb R_+[\Gamma]$ on $X$, where $\Gamma$ is a rational curve of minimal anticanonical degree among curves whose numerical class belongs to $R$;
for any ray $R$ satisfying $L \cdot \Gamma >0$, we define the invariant $\tau_L(R):=-K_X\cdot\Gamma/(L \cdot \Gamma)$ (see \cite[Definition 1.1]{BKLN-LPosRays}) and we call it the {\em $L$-length} of $R$.
This does not require $L$ to be nef, so we can in fact work with any line bundle $L$, {\em i.e.} with any {\em pre-polarized manifold}  $(X,L)$. 
%
In \cite{BKLN-LPosRays} we deal with varieties with extremal rays of $L$-length $> n-2$. 

However, it is possible to consider non integral values of $\tau$. If $L$ is ample, {\em i.e.} for polarized manifolds, we refer to \cite[Chapter 7]{bookBS}, \cite{Fuj-KodEn}, \cite{Zh2-nonint} and \cite{BdT-nonint}. 
In this paper, we investigate the general situation when a pre-polarized variety admits an extremal ray whose $L$-length is not integer, so that it satisfies $n-k<\tau_L(R)<n-k+1$, $n \not= k$, with the technical assumption (cf. \cite[Theorem 2.1]{BdT-nonint}) $n \geq 2k-3$.

The paper is organized as follows:
in Section (\ref{back}) we recall some background material, while
in Section (\ref{contr}) pre-polarized manifolds $(X,L)$ admitting a nef $L$-positive extremal ray of length $\geq n-2$ are described;
building on these descriptions, in Section (\ref{mainnonint}) we classify pairs $(X,L)$ admitting an extremal ray of non-integral $L$-length;
finally, in Section (\ref{mainnefvaluemorphism}) we apply our results to describe rays-positive manifolds (see Definition \ref{rayspos}) and the nefvalue morphism of $(X,L)$ if in addition we assume that $L$ is ample.

\section{Background material}\label{back}

Let $X$ be a smooth projective variety of dimension $n$ defined over the field of complex numbers.
Let $N_1(X)$ be the $\mathbb R$-vector space of $1$-cycles modulo numerical equivalence.
We denote with $\rho_X$ its dimension and we call it the {\em Picard number} of $X$.
Inside $N_1(X)$ we consider the {\em Kleiman--Mori cone} $\overline{\cone}(X)$ of $X$, that is the closure of the cone of the effective $1$-cycles on $X$.
If the negative part of $\overline{\cone}(X)$ (with respect to $K_X$) is not empty, then 
a face in this part of the cone is called an {\em extremal face}, and if it is $1$-dimensional it is called an {\em extremal ray}. 
By the {\em Contraction Theorem}, to any extremal face $\Sigma$ is associated a proper surjective morphism $\Phi\colon X \to Y$ onto a normal variety which exactly contracts all the curves with numerical class in $\Sigma$, whose fibers are connected and such that $-K_X$ is $\Phi$-ample. Such a morphism is usually called a {\em Fano--Mori contraction}, or an {\em extremal contraction}; it is said to be an {\em elementary} contraction when it is associated with an extremal ray.
Moreover, we say that $\Phi$ is {\em of fiber type} if $\dim X > \dim Y$, otherwise we say that it is {\em birational}. In the last case we say that $\Phi$ is {\em divisorial} if it contracts an $(n-1)$-dimensional subvariety of $X$.


An extremal ray is denoted by $R$, its contraction~by $\varphi_R$ and the exceptional~locus of $\varphi_R$ by $E$.
If $\varphi_R$ is of fiber type, we say that $R$ is {\em nef}, otherwise we say that~$R$~is {\em non nef}.
We will write $R$ as $\R_+[\Gamma]$, with~$\Gamma$~a rational curve such that $-K_X \cdot \Gamma=\ell(R)$, where $\ell(R)$ is the {\em length} of $R$ (that is the minimun anticanonical degree of rational curves contracted by $\varphi_R$).

\smallskip

We recall that a smooth complex projective variety $X$ is called a {\em Fano manifold} if its anticanonical bundle $-K_X$ is an ample Cartier divisor.
To a Fano manifold $X$ are associated two invariants, namely the {\em index}, $r_X$, defined as
the largest integer dividing $-K_X$ in the Picard group of $X$,
and the {\em pseudoindex}, $i_X$, defined as
the minimum anticanonical degree of rational curves on $X$. 
Since $X$ is smooth, $\pic(X)$ is torsion free; so the divisor $H$ satisfying $-K_X = r_X H$ is uniquely determined and called the {\em fundamental divisor} of $X$.
It is a classical result that $r_X \leq \dim X+1$,
equality holding if and only if $(X,H)=({\mathbb P}^{\dim X}, {\Ol}_{\mathbb P^{\dim X}}(1))$;
moreover, $r_X={\dim X}$ if and only if $(X,H)=({\mathbb Q}^{\dim X},{\Ol}_{\mathbb Q^{\dim X}}(1))$. 
Finally, we define a {\em del Pezzo manifold} (resp. a {\em Mukai manifold}) as a pair $(X,L)$ where $L$ is an ample line bundle on $X$ such that $-K_X=(\dim X-1)L$ (resp. $-K_X=(\dim X-2)L$).


\smallskip

Throughout the paper we will denote by $\qn$ (resp. $\mathbf Q^n$) the smooth (resp. singular irreducible and reduced) hyperquadric in $\proj^{n+1}$, unless otherwise stated.

\section{Contractions of extremal rays of non-integral $L$-length}\label{contr}
\begin{definition}\cite[cf. Section 1]{BKLN-LPosRays}
Let $X$ be a smooth complex projective variety of dimension $n \geq 3$ and let $L$ be a line bundle on $X$.
An extremal ray $R:=\R_+[\Gamma]$ of $\overline{\cone}(X)$ is said {\em $L$-positive} if $L \cdot \Gamma >0$.
To such a ray we associate the (positive) rational number 
$$\tau_L(R):=\frac{\ell(R)}{L\cdot \Gamma},$$
that we call {\em $L$-length} of $R$. 
\end{definition}

\noindent We will work in the following setup.

\begin{center}
\begin{minipage}[center]{11cm}
\begin{statement}\label{setup} 
Let $X$ be a smooth complex projective variety of dimension $n \geq 3$ and let $L$ be a line bundle on $X$.
Let $R:=\R_+[\Gamma]$ be an $L$-positive extremal ray and denote by $\tau_L(R)$ its $L$-length.
\end{statement}
\end{minipage}
\end{center}\par
\medskip

\noindent We will describe all possible pairs $(X,L)$ under the assumption that $\tau_L(R)\not\in \Z$.
Notice that, since $\tau_L(R) < n+1$, there exists a nonnegative integer $k$ such that
\begin{equation}\label{bounded}
n-k < \tau_L(R) < n-k+1.
\end{equation}
Moreover, being $\tau_L(R)$ a positive number, we have $k \leq n$.\par
\smallskip

\noindent We will make use of these general facts: 
since $L\cdot \Gamma \geq 2$, the length of $R$ is bounded~as
\begin{equation}\label{boundlength}
\ell(R) = \tau_L(R) (L\cdot \Gamma) \geq 2 \tau_L(R) > 2 (n-k),
\end{equation}
hence
\begin{equation}\label{almeno2}
k \geq n - \frac{\ell(R)-1}{2} \Big(\geq \frac{n}{2}\Big), \mbox{ so that }
k \geq 2.
\end{equation}
Furthermore, if $n>k$, the bounds in \eqref{bounded} imply
\begin{equation}\label{boundLGamma}
\ell(R) / (n-k+1) < L \cdot \Gamma < \ell(R) / (n-k).
\end{equation}

\setcounter{equation}{0}

\smallskip

We will make use of the following result.

\begin{lemma}\label{indice}
Let $X$, $L$, $R:=\R_+[\Gamma]$ and $\tau_L(R)$ be as in $(\ref{setup})$.
Assume that $0 \not= n-k < \tau_L(R)<n-k+1$ and that $n \geq 2k-3$.
If $X$ is a Fano manifold with $\rho_X=1$, then $r_X=\ell(R)$ and $L$ is ample.
\end{lemma}

\proof
If $\ell(R)=n+1$, the assertion follows from \cite{CMS-ProjSpace};
so we can assume $\ell(R)\leq n$. Denote by $H$ the fundamental divisor of $X$. 
Since $-K_X=\tau L$, where $\tau:=\tau_L(R)$, we see that $L$ is ample; then $L=mH$ for some positive integer $m$ and we have
\begin{equation}\label{uguali0}
r_X H \cdot \Gamma = -K_X \cdot \Gamma = \tau L \cdot \Gamma = \tau m \ H \cdot \Gamma,
\end{equation}
from which we get $r_X = \tau m$, whence $m \geq 2$.
We have to prove that $H \cdot \Gamma =1$.
\\
Assume, to get a contradiction, that $H \cdot \Gamma\geq 2$. From \eqref{uguali0} we derive
\begin{equation}\label{geq0}
n/2 \geq r_X = \tau m \geq 2 \tau > 2 (n-k),
\end{equation}
which gives $4k \geq 3n+1$.
So we have 
$2 (n+3) \geq 4k \geq 3n+1$, yielding $n \leq 5$.
However, taking into account \eqref{geq0}, this gives a contradiction, as $r_X$ is an integer.
\qed

\setcounter{equation}{0}


\smallskip

There are recurrent situations in our proofs; the easiest case, i.e. when $X$ admits an $L$-positive extremal ray of length $\ell(R)=n+1$, is settled in the following

\begin{proposition}\label{R_n+1}
Let $X$, $L$, $R:=\R_+[\Gamma]$ and $\tau_L(R)$ be as in $(\ref{setup})$.
Assume that $0 \not= n-k < \tau_L(R) < n-k+1$.
If $\ell(R)=n+1$, then one of the following holds:
\begin{itemize}
	\item[(1)] $k\geq 2$ and $(X,L)=(\proj^{2k},\Ol_{\proj^{2k}}(2))$;
	\item[(2)] $k=2$ and $(X,L)=(\proj^3,\Ol_{\proj^3}(3))$;
	\item[(3)] $k=3$ and
\begin{itemize}
	\item[(3-1)] $(X,L)=(\proj^4,\Ol_{\proj^4}(3))$;
	\item[(3-2)] $(X,L)=(\proj^4,\Ol_{\proj^4}(4))$;
\end{itemize}
	\item[(4)] $k=4$ and 
\begin{itemize}
        \item[(4-1)] $(X,L)=(\proj^6,\Ol_{\proj^6}(3))$;
	\item[(4-2)] $(X,L)=(\proj^5,\Ol_{\proj^5}(4))$;
	\item[(4-3)] $(X,L)=(\proj^5,\Ol_{\proj^5}(5))$;
\end{itemize}
	\item[(5)] $k=5$ and $(X,L)=(\proj^7,\Ol_{\proj^7}(3))$;
	\item[(6)] $k=6$ and $(X,L)=(\proj^9,\Ol_{\proj^9}(3))$;
	\item[(7)] $k+1 \leq n \leq 2k-4$ (so $k\geq 5$) and $(X,L)=(\pn,\Ol_\pn(m))$, with $m:=L\cdot\Gamma$.
\end{itemize}
\end{proposition}
\proof
First of all note that, since $\tau_L(R)>0$ and $n-k\not=0$, we have $n>k$.
On the other hand, $n \leq 2k$ by \eqref{boundlength}.
Moreover, the assumption $\ell(R)=n+1$ implies $X=\pn$ by \cite{CMS-ProjSpace}.
%
%
Since, according \eqref{almeno2}, $k \geq 2$,
it is now easy to show the assertion by using the inequalities in \eqref{boundLGamma}.
\qed

\setcounter{equation}{0}


\medskip

If $X$ admits an $L$-positive extremal ray of length $\ell(R)=n$, we have following
\begin{proposition}\label{R_n}
Let $X$, $L$, $R:=\R_+[\Gamma]$ and $\tau_L(R)$ be as in $(\ref{setup})$.
Assume that $0 \not= n-k < \tau_L(R) < n-k+1$ and that $n \geq 2k-3$.
If $\ell(R)=n$, then one of the following holds:
\begin{itemize}
\item[(i)] $\rho_X=1$ and one of the following holds:
\begin{itemize}
	\item[(1)] $k\geq 2$ and
	$(X,L)=(\quadr^{2k-1},\Ol_{\quadr^{2k-1}}(2))$;
	\item[(2)] $k=3$ and $(X,L)=(\quadr^4,\Ol_{\quadr^4}(3))$;
	\item[(3)] $k=4$ and
\begin{itemize}
	\item[(3-1)] $(X,L)=(\quadr^5,\Ol_{\quadr^5}(4))$;
	\item[(3-2)] $(X,L)=(\quadr^5,\Ol_{\quadr^5}(3))$;
\end{itemize}
	\item[(4)] $k=5$ and$(X,L)=(\quadr^7,\Ol_{\quadr^7}(3))$;
\end{itemize}
\item[(ii)] $\rho_X = 2$, $\varphi_R (X)$ is a smooth curve and one of the following holds:
\begin{itemize}
	\item[(5)] $k\geq 2$, and	$(F,L_F)=(\proj^{2k-2},\Ol_{\proj^{2k-2}}(2))$ for any fiber $F$ of $\varphi_R$;
	\item[(6)] $k=3$ and $(F,L_F)=(\proj^3,\Ol_{\proj^3}(3))$ for any fiber $F$ of $\varphi_R$;
	\item[(7)] $k=4$ and
\begin{itemize}
	\item[(7-1)] $(F,L_F)=(\proj^4,\Ol_{\proj^4}(3))$ for any fiber $F$ of $\varphi_R$;
	\item[(7-2)] $(F,L_F)=(\proj^4,\Ol_{\proj^4}(4))$ for any fiber $F$ of $\varphi_R$;
\end{itemize}
	\item[(8)] $k=5$, and	$(F,L_F)=(\proj^6,\Ol_{\proj^6}(3))$ for any fiber $F$ of $\varphi_R$.
\end{itemize}
\end{itemize}
\end{proposition}

\proof
The assumption $\ell(R)=n$ together with equation \eqref{boundlength} implies $n \leq 2k-1$; moreover, \cite[Proposition 2.4]{Wi-Length} gives the following possibilities:
\begin{itemize}
	\item[(i)] $X$ is a Fano manifold with $\rho_X=1$;
	\item[(ii)] $\rho_X=2$ and $\varphi_R \colon X \to Y$ is a morphism onto a smooth curve $Y$ whose general fiber $F$ is a smooth variety of dimension $n-1$ and Picard number $\rho_F=1$ admitting an extremal ray $R_{(F)}:=\R_+[C]$ of length $\dim F+1=n$.
\end{itemize}

Assume first that $X$ is as in case (i)
so that,by Lemma (\ref{indice}), $X$ is a Fano manifold with $r_X=\dim X$; it follows that $(X,H)=(\qn,\Ol_\qn(1))$.
%
%
Recall that, by \eqref{almeno2}, $k \geq 2$.
If $k=2$, then $n=3$; hence taking into account the inequalities in \eqref{boundLGamma} we get case (1) (with $k=2$) of the statement;
if $k=3$, then $4 \leq n \leq 5$, so, using \eqref{boundLGamma}, we derive the cases (2) and (1) (with $k=3$);
as to $k \geq 4$, recall that $2k-3 \leq n \leq 2k-1$,
so, by \eqref{boundLGamma}, we obtain cases (1) (with $k\geq 4$),
(3) and~(4).

Assume now that $X$ is as in case (ii) and consider $F$, the general fiber of $\varphi_R$.
Since $K_X+\tau_L(R) L$ is trivial on $F$, we have $\tau_{(F)}:=\tau_{L_F}(R_{(F)})=\tau_L(R)$, so
$\dim F - k_{(F)} < \tau_{(F)} < \dim F - k_{(F)} + 1,$
where $k_{(F)}:=k-1$.
If $k=2$, then $n=3$; so $L_F \cdot C=2$, in view of the bounds in \eqref{boundLGamma}. It follows that the restriction of $2K_X+3L$ to $F$ is trivial, hence it is immediate to derive $(F,L_F)=(\pd,\Ol_\pd(2))$. Moreover, all fibers of $\varphi_R$ are irreducible and reduced, so by semicontinuity we find that, for any fiber $G$ of $\varphi_R$, $0 \leq \Delta(G,L_G) \leq \Delta(F,L_F) = 0$, whence $(G,L_G)=(\pd,\Ol_\pd(2))$ by \cite[Theorems 5.10 and 5.15]{bookFuj}. This leads to case (5) (with $k=2$) of the statement.
As to $k \geq 3$, we apply Proposition (\ref{R_n+1}) to $F$.
Therefore $(F,L_F)=(\proj^{n-1}, \Ol_{\proj^{n-1}}(a))$, where $a$ is known.
If $(F,L_F)\not=(\proj^4, \Ol_{\proj^4}(3))$, the line bundle $N:=-K_X-(n-k)L$ is $\varphi_R$-ample and $N\cdot \Gamma=1$.
Since, for some~ample line bundle $A$ on $Y$, the line bundle $M:=N+\varphi_R^\ast A$ is ample on $X$ and $K_X+nM$ is trivial on $R$, we derive $(F, M_F)=(\proj^{n-1}, \Ol_{\proj^{n-1}}(1))$. 
Since $\varphi_R$ is equidimensional,  $(F, M_F)=(\proj^{n-1},$ $ \Ol_{\proj^{n-1}}(1))$ for any fiber of $\varphi_R$ by \cite[Lemma 2.12]{Fuj-Polarized}. It is now immediate to see~that $(F,L_F)=(\proj^{n-1}, \Ol_{\proj^{n-1}}(a))$ for any fiber $F$ of $\varphi_R$.
If $(F,L_F)=(\proj^4, \Ol_{\proj^4}(3))$,~we consider the line bundle $N:=K_X+2L$ and we conclude with the same~argument.
%
\qed

\medskip

If $X$ admits an $L$-positive extremal ray of length $\ell(R)=n-1$, we have following
\begin{proposition}\label{R_n-1}
Let $X$, $L$, $R:=\R_+[\Gamma]$ and $\tau_L(R)$ be as in $(\ref{setup})$.
Assume that $0 \not= n-k < \tau_L(R) < n-k+1$ and that $n \geq 2k-3$.
If $R$ is nef and $\ell(R)=n-1$, then one of the following holds:
\begin{itemize}
\item[(i)] $\rho_X=1$ and one of the following holds:
\begin{itemize}
	\item[(1)] $k\geq 3$, $(X,H)$ is a del Pezzo $(2k-2)$-fold and $L=2H$;
	\item[(2)] $k=4$, $(X,H)$ is a del Pezzo $5$-fold and $L=3H$;
\end{itemize}
\item[(ii)] $\rho_X \geq 2$, $\varphi_R (X)$ is a smooth variety and one of the following holds:
\begin{itemize}
 \item[(3)] $k\geq 3$, $n=2k-2$ and
\begin{itemize}
	\item[(3-1)] $(F,L_F)=(\proj^{2k-4},\Ol_{\proj^{2k-4}}(2))$ for any fiber $F$;
	\item[(3-2)] $(F,L_F)=(\quadr^{2k-3},\Ol_{\quadr^{2k-3}}(2))$ for all smooth fibers $F$, $(G,L_G)=({\mathbf Q}^{2k-3},\Ol_{{\mathbf Q}^{2k-3}}(2))$ for the singular fibers $G$ (if any);
\end{itemize}
 \item[(4)] $k=4$, $n=5$ and
\begin{itemize}
	\item[(4-1)] $(F,L_F)=(\pt,\Ol_{\pt}(3))$ for any fiber $F$;
	\item[(4-2)] $(F,L_F)=(\quadr^4,\Ol_{\quadr^4}(3))$ for all smooth fibers $F$, $(G,L_G)=({\mathbf Q}^4,$ $\Ol_{{\mathbf Q}^4}(3))$ for the singular fibers $G$ (if any);
\end{itemize}
\end{itemize}
\end{itemize}
\end{proposition}
\proof
By assumption $\ell(R)=n-1$, so $n \leq 2k-2$ from equation \eqref{boundlength}.
Note that $k \geq 3$, since $n \geq 3$.
Moreover, by \eqref{boundLGamma}, if $n=2k-3$ then $L\cdot\Gamma=3$ and $k=4$, while if $n=2k-2$ then $L\cdot\Gamma=2$.

Assume first that $\rho_X=1$; 
then, according to Lemma (\ref{indice}), $X$ is a Fano manifold with $r_X=n-1$; so $(X,H)$ is a del Pezzo manifold.
%
Since $L=3H$ if $n=2k-3$, and $L=2H$ if $n=2k-2$, it is straightforward to get case (i) of the statement.

We can thus assume that $\rho_X \geq 2$. By \cite[Theorem 1.1]{Wi-OnContr} the target of the contraction $\varphi_R\colon X \to Y$ is a variety of dimension $1$ or $2$; hence $Y$ is smooth and $\varphi_R$ is equidimensional by \cite[Cf. Proposition 1.4.1]{ABW-VBAdj}.
Since $K_X+\tau_L(R) L$ is trivial on general fibers $F$ of $\varphi_R$, these fibers are Fano manifolds of pseudoindex $i_F \geq \dim F$; so $\rho_F=1$ by \cite[Theorem A]{Wi-Mukai}.
Moreover, $F$ admits an extremal ray $R_{(F)}$ such that $\ell(R_{(F)})=n-1$, hence $\tau_{(F)}:=\tau_{L_F}(R_{(F)})=\tau_L(R)$, so
$\dim F - k_{(F)} < \tau_{(F)} < \dim F - k_{(F)} + 1,$
where $k_{(F)}:=k-\dim Y$.
Recall that $n=2k-3$ or $2k-2$.
\\
If $\dim Y=1$, then $\dim F=\ell(R_{(F)})$ and $k_{(F)}:=k-1 (\geq 2)$. 
So Proposition (\ref{R_n}) applies and $(F,L_F)$ is either
$(\quadr^4,\Ol_{\quadr^4}(3))$ if $k=4$, or $(\quadr^{2k-3},$ $\Ol_{\quadr^{2k-3}}(2))$ if $k\geq 3$. 
In the first case, the line bundle $N:=-K_X-L$ is $\varphi_R$-ample and $N\cdot \Gamma=1$.
Since $\varphi_R$ is equidimensional and, for some ample line bundle $A$ on $Y$, $M:=N+\varphi_R^\ast A$ is ample and $K_X+(n-1)M$ is nef and it is trivial only on $R$, by \cite[Theorem B]{ABW-TwoThms} we have that ${\mathcal E}:=$ $\varphi_{R \ast} M$ is a locally free sheaf of rank $n+1$ and $X$ embeds into $\proj_Y({\mathcal E})$ as a divisor of relative degree $2$. Then we derive case (4-2) of the statement.
In the last case, the same argument applied to $N:=K_X+(k-1)L$ give case~(3-2).
\\
If $\dim Y=2$, then $\dim F=\ell(R_{(F)})-1$ and $k_{(F)}:=k-2 (\geq 1)$.
If $k=3$, then $n=4$ and, since $2K_X+3L$ is trivial on $F$, it is immediate to derive $(F,L_F)=(\pd,\Ol_\pd(2))$;
if $k \geq 4$, then Proposition (\ref{R_n+1}) applies, so $(F,L_F)$ can be
$(\pt,\Ol_\pt(3))$ if $k=4$ and $(\proj^{2k-4},\Ol_{\proj^{2k-4}}(2))$ for $k\geq 4$.
If the first two cases, the line bundle $N:=-K_X-L$ is $\varphi_R$-ample and $N\cdot \Gamma=1$.
Since, for some ample line bundle $A$ on $Y$, the line bundle $M:=N+\varphi_R^\ast A$ is ample and $K_X+(n-1)M$ is trivial on $R$, we derive that $(F, M_F)=(\proj^{n-2}, \Ol_{\proj^{n-2}}(1))$. 
Since $\varphi_R$ is equidimensional,  $(F, M_F)=(\proj^{n-2}, \Ol_{\proj^{n-2}}(1))$ for any fiber of $\varphi_R$ by \cite[Lemma 2.12]{Fuj-Polarized}. It is immediate to see that, for any fiber $F$ of $\varphi_R$, $(F,L_F)=(\proj^{n-1}, \Ol_{\proj^{n-2}}(a))$, with $a=2$ and $3$, resp.. 
In the last case, we can apply the same argument to $N:=K_X+(k-1)L$.
So we get cases (4-1) and (3-1).
\qed

\medskip

If $X$ admits an $L$-positive extremal ray of length $\ell(R)=n-2$, we have following
\begin{proposition}\label{R_n-2}
Let $X$, $L$, $R:=\R_+[\Gamma]$ and $\tau_L(R)$ be as in $(\ref{setup})$.
Assume that $0 \not= n-k < \tau_L(R) < n-k+1$ and that $n \geq 2k-3$.
If $R$ is nef and $\ell(R)=n-2$, then one of the following holds:
\begin{itemize}
\item[(i)] $\rho_X=1$,
\begin{itemize}
	\item[(1)] $k\geq 4$, $(X,H)$ is a Mukai $(2k-3)$-fold and $L=2H$; 
\end{itemize}
\item[(ii)] $\rho_X \geq 2$ and one of the following holds:
\begin{itemize}
 \item[(2)] $k=4$, $\varphi_R(X)$ is a smooth curve and $(F,L_F)=(\pd\times\pd, \Ol_{\pd\times\pd}(2,2))$ for the general fiber $F$ of $\varphi_R$;
 \item[(3)] $k\geq 4$, 
\begin{itemize}
	\item[(3-1)] $\varphi_R(X)$ is a smooth curve, $(F,H_{(F)})$ is del Pezzo $(2k-4)$-fold with $\rho_F=1$ and $L_F=2H_{(F)}$ for the general fiber $F$ of $\varphi_R$;
	\item[(3-2)] $\varphi_R(X)$ is a smooth surface, $(F,L_F)=(\quadr^{2k-5},\Ol_{\quadr^{2k-5}}(2))$ for all smooth fibers $F$ of $\varphi_R$ and $(G,L_G)=(\mathbf Q^{2k-5}, \Ol_{\mathbf Q^{2k-5}}(2))$ for singular fibers $G$ (if any);
	\item[(3-3)] $\varphi_R(X)$ is a $3$-fold with at most isolated rational Gorenstein singularities, $(F,L_F)=(\proj^{2k-6},\Ol_{\proj^{2k-6}}(2))$ for any fiber $F$ over the smooth locus of $\varphi_R(X)$ and $\dim G=2k-5$ for all fibers $G$ over the singular locus of $\varphi_R$ (if any);
\end{itemize}
\end{itemize}
\end{itemize}
\end{proposition}

\proof
Notice that our assumptions together with condition (\ref{boundlength}) imply $n=2k-3$, so that $k \geq 4$.
Moreover, by (\ref{boundLGamma}) we have $L\cdot\Gamma=2$.

Assume first that $\rho_X=1$;
then, according to Lemma (\ref{indice}), $X$ is a Fano manifold with $r_X=n-2$, so $(X,H)$ is a Mukai manifold and
we get case (i) of the statement.

We can thus assume that $\rho_X \geq 2$. By \cite[Theorem 1.1]{Wi-OnContr} the target of the contraction $\varphi_R\colon X \to Y$ is a variety of dimension $1$, or $2$, or~$3$.
If $\dim Y=1$ or $2$, then $Y$ is smooth and $\varphi_R$ is equidimensional by \cite[Cf. Proposition 1.4.1]{ABW-VBAdj}, while, if $\dim Y=3$, then it is well-known that $Y$ has rational, Gorenstein singularities (cf. \cite[Corollary 7.4]{Ko-Dualizing}).
Moreover, they are also isolated: take a general divisor $D\in |\varphi_R^\ast\Ol_Y(1)|$; since $\varphi_R^\ast\Ol_Y(1)$ is base point free and $X$ is smooth,~by Bertini's theorem also $D$ is smooth; then $\varphi_R$ restricted to $D$ is an elementary contraction onto a surface $S$, which is smooth by \cite[Cf. Proposition 1.4.1]{ABW-VBAdj}; so $Y$ has isolated singularities.
\\
If $k=4$, then $n=5$ and $\ell(R)=3$. Moreover, the general fiber $F$ is a Fano manifold such that $2 \leq \dim F \leq 4$ and $i_F=3$.
Therefore, unless $F=\pd\times\pd$, $\rho_F=1$ (see \cite[Theorem 3]{NO-RCBound}). 
In this last case, let $H_{(F)}$ be the fundamental divisor of $F$. 
Since $K_X+\frac{3}{2} L$ is trivial on $F$, then $(F,H_{(F)})$ is one of the following:
a del Pezzo $4$-fold, $(\qt,\Ol_\qt(1))$ (by \cite[Corollary 2.6]{Wi-Length}) or $(\pd,\Ol_\pd(1))$. Therefore we get cases (2) and (3) (with $k=4$) of the statement.
\\
As to $k \geq 5$, since the general fiber $F$ of $\varphi_R$ admits an extremal ray $R_{(F)}$ of length $n-2$, we have 
$\dim F - k_{(F)} < \tau_{(F)} < \dim F - k_{(F)} + 1,$
where $\tau_{(F)}:=\tau_{L_F}(R_{(F)})=\tau_L(R)$ and $k_{(F)}:=k-\dim Y$.
Since $K_X+\tau_L(R) L$ is trivial on $F$, such $F$ is a Fano manifold with $i_F\geq \dim F-1$; so, recalling that $\dim F \geq n-3\geq 4$, by using \cite[Theorem A]{Wi-Mukai} and \cite[Theorem 1]{Occ-Products}
, we derive $\rho_F=1$ unless $F=\pd\times\pd$. 
However, in this last case, $k=5$ and $n=7$, so $F$ cannot have extremal rays of length $n-2=5$.
Then $\rho_F=1$.
\\
If $\dim Y=1$, then, for a suitable ample line bundle $H_{(F)}$ on $F$, $(F,H_{(F)})$ is 
a del Pezzo $(2k-4)$-fold and $L_F=2H_{(F)}$ by Proposition (\ref{R_n-1}); so we get case (3-1) (with $k\geq 5$) of the statement.
\\
If $\dim Y=2$, then $(F,L_F)=(\quadr^{2k-5},\Ol_{\quadr^{2k-5}}(2))$ by Proposition (\ref{R_n}), while
if $\dim Y=3$, then $(F,L_F)=(\proj^{2k-6},\Ol_{\proj^{2k-6}}(2))$ by Proposition (\ref{R_n+1}).
To complete the proof, note that in both cases the line bundle $N:=K_X+(k-2)L$ is $\varphi_R$-ample and $N\cdot \Gamma=1$.
Moreover, for some ample line bundle $A$ on $Y$, the line bundle $M:=N+\varphi_R^\ast A$ is ample and $K_X+(n-2)M$ is nef and it is trivial only on $R$.
\\
If $\dim Y=2$, since $\varphi_R$ is equidimensional, by \cite[Theorem B]{ABW-TwoThms} we have that ${\mathcal E}:=\varphi_{R \ast} M$ is a locally free sheaf of rank $n$ and $X$ embeds into $\proj_Y({\mathcal E})$ as a divisor of relative degree $2$, so we get case (3-2) (with $k\geq 5$). 
\\
If $\dim Y=3$, following the proof of \cite[Lemma 2.12]{Fuj-Polarized} we see that $\varphi_R(F)$ is a smooth point if $\dim F=n-3$ since $\varphi_R$ cannot have divisorial fibers by \cite[Lemma 1.2]{ABW-OnContr} and the argument is local.
Therefore the dimension of any fiber over the singular locus of $Y$ is $n-2$. 
Moreover, there does not exist any fiber over the smooth locus of $Y$ whose dimension in $n-2$, since otherwise $n$ would be equal to $4$ by \cite[Theorem 4.1]{AW-NonVanish}, contradicting $n > k$.
Therefore $\varphi_R$ is equidimensional over the smooth locus of $Y$, so that $(X,M)$ is a scroll here (cf. the proof of \cite[Lemma 2.12]{Fuj-Polarized}).
It is now immediate to see that, for any fiber $F$ of $\varphi_R$ over the smooth locus, $(F,L_F)=(\proj^{n-1},$ $\Ol_{\proj^{n-2}}(2))$.
Then we get case (3-3) (with $k\geq 5$) of the~statement.
\qed

\setcounter{equation}{0}

\section{Varieties with an extremal ray of non-integral $L$-length}\label{mainnonint}

In this section we apply the results of the previous section to classify pairs $(X,L)$, where $X$, $L$ are as in $(\ref{setup})$ and \eqref{bounded} holds.
In particular, according to condition \eqref{almeno2}, $k\geq 2$. Case $k=2$ is already settled in \cite[Proposition 1.4]{BKLN-LPosRays}, so we confine to $k \geq 3$. 
The following proposition takes care of the case $k=3$.

\begin{proposition}\label{noninteger3} 
Let $X$, $L$, $R:=\R_+[\Gamma]$ and $\tau_L(R)$ be as in $(\ref{setup})$.
Assume that $0 \not= n-3 < \tau_L(R) < n-2$.
Then one of the following holds:
\begin{itemize}
\item[(i)] $\rho_X=1$ and one of the following holds:
\begin{itemize}
\item[(1)] $n=6$, 
	$\tau_L(R)=\frac{7}{2}$ and $(X,L)=(\proj^6,\Ol_{\proj^6}(2))$;
\item[(2)] $n=5$, 
	$\tau_L(R)=\frac{5}{2}$ and $(X,L)=(\quadr^5,\Ol_{\quadr^5}(2))$;
\item[(3)] $n=4$, 
\begin{itemize}
	\item[(3-1)] $\tau_L(R)=\frac{5}{4}$ and $(X,L)=(\proj^4,\Ol_{\proj^4}(4))$;
	\item[(3-2)] $\tau_L(R)=\frac{5}{3}$ and $(X,L)=(\proj^4,\Ol_{\proj^4}(3))$;
	\item[(3-3)] $\tau_L(R)=\frac{4}{3}$ and $(X,L)=(\quadr^4,\Ol_{\proj^4}(3))$;
	\item[(3-4)] $\tau_L(R)=\frac{3}{2}$, $(X,H)$ is a del Pezzo manifold and $L=2H$;
\end{itemize}
\end{itemize}
\item[(ii)] $\rho_X \geq 2$, $\varphi_R(X)$ is a smooth variety and one of the following holds:
\begin{itemize}
\item[(4)] $n=5$, 
  $\tau_L(R)=\frac{5}{2}$ and $(F,L_F)=(\proj^4,\Ol_{\proj^4}(2))$ for any fiber $F$ of~$\varphi_R$;
\item[(5)] $n=4$, 
\begin{itemize}
\item[(5-1)] $\tau_L(R)=\frac{4}{3}$ and $(F,L_F)=(\pt,\Ol_{\pt}(3))$ for any fiber $F$ of $\varphi_R$;
\item[(5-2)] $\tau_L(R)=\frac{3}{2}$ and $(F,L_F)=(\quadr^3,\Ol_{\quadr^3}(2))$ for all smooth fiber $F$ of $\varphi_R$ and $(G,L_G)=({\mathbf Q}^3,\Ol_{{\mathbf Q}^3}(2))$ for singular fibers $G$ (if any);
\item[(5-3)] $\tau_L(R)=\frac{3}{2}$ and $(F,L_F)=(\pd,\Ol_{\pd}(2))$ for any fiber $F$ of $\varphi_R$;
\item[(5-4)] $\tau_L(R)=\frac{3}{2}$ and $\varphi_R$ is the blow-up of $Y$ at one point.
\end{itemize}
\end{itemize}
\end{itemize}
\end{proposition}

\proof
Assume first that $R$ is a non nef extremal ray. Then the bound in \eqref{boundlength} combined with \cite[Theorem 1.1]{Wi-OnContr} yields $n=4$ and $\ell(R)=3$. Therefore, by \cite[Theorem 1.1]{AO-SpeRays}, $\varphi_R\colon X \to Y$ is the blow-up of a smooth $4$-fold $Y$ at one point.
By computing $L \cdot \Gamma$ with \eqref{boundLGamma}, we get case (5-4) in the statement.

Suppose now that $R$ is nef. Then the bound in \eqref{boundlength} implies $n\leq 6$.
If $n=6$, then $\ell(R)=7$, so we get case (1) of the statement by Proposition (\ref{R_n+1}).
If $n=5$, then $6\geq \ell(R)\geq 5$, so we get cases (2) and (4) by Propositions (\ref{R_n+1}) and (\ref{R_n}).
If $n=4$, then $5\geq \ell(R)\geq 3$, so we get cases (3) and (5-1)--(5-3) by Propositions (\ref{R_n+1}), (\ref{R_n}) and (\ref{R_n-1}).
\qed

\smallskip

Now, we can assume $k \geq 4$. We obtain the following

\begin{proposition}\label{nonintegerk} 
Let $X$, $L$, $R:=\R_+[\Gamma]$ and $\tau_L(R)$ be as in $(\ref{setup})$.
Assume that $0 \not= n-k < \tau_L(R) < n-k+1$, $k \geq 4$ and $n \geq 2k-3$.
Then one of the following holds:
\begin{itemize}
\item[(i)] $\rho_X=1$ and one of the following holds:
\begin{itemize}
\item[(1)] $n=2k$, 
	$\tau_L(R)=\frac{n+1}{2}$ and $(X,L)=(\pn,\Ol_\pn(2))$;
\item[(2)] $n=2k-1$, 
	$\tau_L(R)=\frac{n}{2}$ and $(X,L)=(\qn,\Ol_\qn(2))$;
\item[(3)] $n=2k-2$, 
\begin{itemize}
	\item[(3-1)] $\tau_L(R)=\frac{7}{3}$ and $(X,L)=(\proj^6,\Ol_{\proj^6}(3))$;
	\item[(3-2)] $\tau_L(R)=\frac{n-1}{2}$, $(X,H)$ is a del Pezzo manifold and $L=2H$;
\end{itemize}
\item[(4)] $n=2k-3$, 
\begin{itemize}
	\item[(4-1)] $\tau_L(R)=\frac{10}{3}$ and $(X,L)=(\proj^9,\Ol_{\proj^9}(3))$;
	\item[(4-2)] $\tau_L(R)=\frac{8}{3}$ and $(X,L)=(\proj^7,\Ol_{\proj^7}(3))$;
	\item[(4-3)] $\tau_L(R)=\frac{6}{5}$ and $(X,L)=(\proj^5,\Ol_{\proj^5}(5))$;
	\item[(4-4)] $\tau_L(R)=\frac{3}{2}$ and $(X,L)=(\proj^5,\Ol_{\proj^5}(4))$;
	\item[(4-5)] $\tau_L(R)=\frac{7}{3}$ and $(X,L)=(\quadr^7,\Ol_{\quadr^7}(3))$;
	\item[(4-6)] $\tau_L(R)=\frac{5}{4}$ and $(X,L)=(\quadr^5,\Ol_{\quadr^5}(4))$;
	\item[(4-7)] $\tau_L(R)=\frac{5}{3}$ and $(X,L)=(\quadr^5,\Ol_{\quadr^5}(3))$;
	\item[(4-8)] $\tau_L(R)=\frac{4}{3}$, $(X,H)$ is a del Pezzo $5$-fold and $L=3H$;
	\item[(4-9)] $\tau_L(R)=\frac{n-2}{2}$, $(X,H)$ is a Mukai manifold and $L=2H$;
\end{itemize}
\end{itemize}
\item[(ii)] $\rho_X \geq 2$ and one of the following holds:
\begin{itemize}
\item[(5)] $n=2k-1$, 
 $\tau_L(R)=\frac{n}{2}$, $\varphi_R(X)$ is a smooth curve and $(F,L_F)=(\proj^{n-1},$ $\Ol_{\proj^{n-1}}(2))$ for any fiber $F$ of $\varphi_R$;
\item[(6)] $n=2k-2$, 
\begin{itemize}
\item[(6-1)] $\tau_L(R)=\frac{n-1}{2}$, $\varphi_R(X)$ is a smooth surface and $(F,L_F)=(\proj^{n-2},$ $\Ol_{\proj^{n-2}}(2))$ for any fiber $F$ of~$\varphi_R$;
\item[(6-2)] $\tau_L(R)=\frac{n-1}{2}$, $\varphi_R(X)$ is a smooth curve, $(F,L_F)=(\quadr^{n-1},$ $\Ol_{\quadr^{n-1}}(2))$ for all smooth fibers $F$ of $\varphi_R$ and $(G,L_G)=({\mathbf Q}^{n-1},$ $\Ol_{{\mathbf Q}^{n-1}}(2))$ for~the singular fibers $G$ (if any);
\item[(6-3)] $\tau_L(R)=\frac{n-1}{2}$ and $\varphi_R$ is the blow-up of $Y$ at one point;
\end{itemize}
\item[(7)] $n=2k-3$, 
\begin{itemize}
\item[(7-1)] $\tau_L(R)=\frac{7}{3}$, $\varphi_R(X)$ is a smooth curve and $(F,L_F)=(\proj^6,\Ol_{\proj^4}(3))$ for any fiber $F$ of $\varphi_R$;
\item[(7-2)] $\tau_L(R)=\frac{5}{4}$, $\varphi_R(X)$ is a smooth curve and $(F,L_F)=(\proj^4,\Ol_{\proj^4}(4))$ for any fiber $F$ of $\varphi_R$;
\item[(7-3)] $\tau_L(R)=\frac{5}{3}$, $\varphi_R(X)$ is a smooth curve and $(F,L_F)=(\proj^4,\Ol_{\proj^4}(3))$ for any fiber $F$ of $\varphi_R$;
\item[(7-4)] $\tau_L(R)=\frac{4}{3}$, $\varphi_R(X)$ is a smooth surface and $(F,L_F)=(\pt,$ $\Ol_\pt(3))$ for any fiber $F$ of $\varphi_R$;
\item[(7-5)] $\tau_L(R)=\frac{4}{3}$, $\varphi_R(X)$ is a smooth curve and $(F,L_F)=(\quadr^4,\Ol_{\quadr^4}(3))$ for all smooth fibers $F$ of $\varphi_R$ and $(G,L_G)=({\mathbf Q}^4, \Ol_{{\mathbf Q}^4}(3))$ for~the singular fibers $G$ (if any);
\item[(7-6)] $\tau_L(R)=\frac{3}{2}$, $\varphi_R(X)$ is a smooth curve and $(F,L_F)=(\pd\times\pd, \Ol_{\pd\times\pd}(2,2))$ for the general fiber $F$ of $\varphi_R$;
\item[(7-7)] $\tau_L(R)=\frac{n-2}{2}$, $\varphi_R(X)$ is a smooth curve, $(F,H_{(F)})$ is del Pezzo $(n-1)$-fold with $\rho_F=1$ and $L_F=2H_{(F)}$ for the general fiber $F$ of $\varphi_R$;
\item[(7-8)] $\tau_L(R)=\frac{n-2}{2}$, $\varphi_R(X)$ is a smooth surface, $(F,L_F)=(\quadr^{n-2},$ $\Ol_{\quadr^{n-2}}(2))$ for all smooth fibers $F$ of $\varphi_R$ and $(G,L_G)=(\mathbf Q^{n-2},$ $\Ol_{\mathbf Q^{n-2}}(2))$ for singular fibers $G$ (if any);
\item[(7-9)] $\tau_L(R)=\frac{n-2}{2}$, $\varphi_R(X)$ is a $3$-fold with at most isolated rational Gorenstein singularities, $(F,L_F)=(\proj^{n-3},\Ol_{\proj^{n-3}}(2))$ for any fiber $F$ over the smooth locus of $\varphi_R(X)$ and $\dim G=n-2$ for all fibers $G$ over the singular locus of $\varphi_R$ (if any);
\item[(7-10)] $\tau_L(R)=\frac{4}{3}$ and $\varphi_R$ is the blow-up of a smooth $5$-fold $Y$ at one point;
\item[(7-11)] $\tau_L(R)=\frac{n-2}{2}$ and $\varphi_R$ is the blow-up of a smooth variety $Y$ along a smooth curve;
\item[(7-12)] $\tau_L(R)=\frac{n-2}{2}$, $\varphi_R(E)$ is a point and $(E,-E_E)$ is either $(\proj^{n-1},$ $\Ol_{\proj^{n-1}}(2))$, or
$(\quadr^{n-1},\Ol_{\quadr^{n-1}}(1))$, where $\quadr^{n-1}$ is a possibly singular hyperquadric.
\end{itemize}
\end{itemize}
\end{itemize}
\end{proposition}

\proof
Assume first that $R$ is a non nef extremal ray. Since $\ell(R)\leq n-1$, the bound in \eqref{boundlength} yields $n=2k-2$ or $2k-3$.
In the former case, by \cite[Theorem 1.1]{AO-SpeRays}, we get case (6-3) of the statement, while in the latter, taking into account \eqref{boundLGamma}, we get cases (7-10)--(7-12) by \cite[Theorems 1.1 and 5.2]{AO-SpeRays}.

Assume now that $R$ is nef. Then the bound in \eqref{boundlength} gives $2k-3 \leq n \leq 2k$.
If $n=2k$, then $\ell(R)=n+1$, so we get case (1) of the statement by Proposition~(\ref{R_n+1}).
If $n=2k-1$, then, recalling \eqref{boundLGamma}, $\ell(R)=n$, so we get cases (2) and (5) by Proposition~(\ref{R_n}).
If $n=2k-2$, then, recalling \eqref{boundLGamma}, either $\ell(R)= n+1$, 
or $\ell(R)= n-1$;
in the former case we get case (3-1) by Proposition (\ref{R_n+1}), while in the latter we get cases (3-2), (6-1) and (6-2) by Proposition (\ref{R_n-1}).
If $n=2k-3$, then $n-2 \leq \ell(R) \leq n+1$;
if $\ell(R)=n+1$,~
we get cases (4-1)--(4-4) by Proposition (\ref{R_n+1});
if $\ell(R)=n$, 
we get cases (4-5)--(4-7)~and (7-1)--(7-3) by Proposition (\ref{R_n}); 
if $\ell(R)=n-1$, 
we get cases (4-8), (7-4) and (7-5) by Proposition (\ref{R_n-1});
if $\ell(R)=n-2$, 
we get cases (4-9) and (7-6)--(7-9) by Proposition (\ref{R_n-2}). 
\qed

\section{Application to semi-polarized varieties}\label{mainnefvaluemorphism}

In this section we will use the previous results to describe {\em rays-positive manifolds} (see Definition (\ref{rayspos})) and to describe the nefvalue morphism of polarized manifolds.

\noindent We first recall the notion of {\em rays-positive manifold} introduced in \cite[Section~2]{BKLN-LPosRays}.\par

\smallskip

Let $(X,L)$ be a {\em semi-polarized manifold}, {\em i.e.} a pair consisting of a smooth complex projective variety $X$ of dimension $n\geq 3$ and a nef line bundle $L$ on $X$.
Assume that $K_X$ is not nef, or, equivalently, that $X$ is not minimal in the sense of the {\em Minimal Model Program}.
We can define the numerical invariant given by
\begin{equation}\label{defsigma}
\sigma := \sigma (X,L) = \mbox{sup} \{t \in \R: t K_X+ L \mbox{ is nef}\}.
\end{equation}
By the {\em Kawamata's Rationality Theorem}, $\sigma$ is a (non-negative) rational number and there exists an extremal ray $R$ in $\overline{\cone}(X)$ such that 
$(\sigma K_X+ L) \cdot R=0$.
\\
Clearly $\sigma >0$ if $L$ is ample. We refer to \cite[Section 2]{BKLN-LPosRays} for examples with $\sigma=0$.

\medskip

For our purpose, we will use the following slight modification of \cite[Definition~2.3]{BKLN-LPosRays}:
\begin{definition}(cf. \cite[Definition 2.3 and Lemma 2.5]{BKLN-LPosRays})\label{rayspos}
Let $(X,L)$ be a semi-polarized manifold such that $K_X$ is not nef.
We say that $(X,L)$ is {\em rays-positive} and that $L$ is {\em rays-positive} if $\sigma(X,L) >0$.
\end{definition}

\begin{remark}
If $L$ is ample or numerically positive, then $(X,L)$ is rays-positive. The converses are not true and we refer to \cite[Section 2]{BKLN-LPosRays} for further discussion and examples.
\end{remark}
 
Notice that, for any extremal ray $R$ orthogonal to $\sigma(X,L)K_X+L$, we have $\tau_L(R)=\frac{1}{\sigma(X,L)}$.
Moreover, since $n \geq 3$ and we are interested in non-integral values of $\tau_L(R)$, in view of \cite[Proposition 1.2]{BKLN-LPosRays}, we have $\frac{1}{\sigma(X,L)}< n-1$.

\medskip

The following proposition deals with rays-positive manifolds such that $n-2 < \frac{1}{\sigma(X,L)} < n-1$.

\begin{proposition}
Let $(X,L)$ be a rays-positive manifold of dimension $n$. Assume that $n-2 < \frac{1}{\sigma(X,L)} < n-1$.
Then $(X,L)$ admits exactly one extremal ray satisfying $(\sigma(X,L) K_X+ L) \cdot R=0$ and is described in \cite[Proposition 1.4]{BKLN-LPosRays}.
\end{proposition}

\proof
Let $R$ be an extremal ray satisfying $(\sigma(X,L) K_X+ L) \cdot R=0$; then we are in the assumptions of \cite[Proposition 1.4]{BKLN-LPosRays}.
Since fibers of different extremal rays can meet only at points, it is immediate to see that $X$ admits only one of such rays.
\qed

\medskip

In the following proposition we take care of rays-positive manifolds such that $0 \not= n-3 < \frac{1}{\sigma(X,L)} < n-2$.

\begin{proposition}
Let $(X,L)$ be a rays-positive manifold of dimension $n$. Assume that $0 \not= n-3 < \frac{1}{\sigma(X,L)} < n-2$.
Then $(X,L)$ admits exactly one extremal ray satisfying $(\sigma(X,L) K_X+ L) \cdot R=0$ and it is described in Proposition (\ref{noninteger3}), unless $\sigma(X,L)=\frac{2}{3}$ and one of the following holds:
\begin{itemize}
	\item[(1)] $X=\pd\times\pd$ and $L$ restricts as $\Ol_\pd(2)$ on any fiber of each projection;
	\item[(2)] $X$ is the blow-up of the smooth variety $Y$ at a finite number of points and $(E_i,L_{E_i})=(\pt,\Ol_\pt(2))$
	 for all exceptional divisors $E_i$.
\end{itemize}
\end{proposition}

\proof
Let $\{R_i:=\mathbb R_+ [\Gamma_i]\}_{i\in I}$ be the set of all the extremal rays in $\overline{\cone}(X)$ satisfying $(\sigma(X,L) K_X+ L) \cdot R_i=0$ and denote by $\tau:=\tau_L(R_i)=\frac{1}{\sigma(X,L)}$.
Then each $R_i$ satisfies the assumptions of Proposition (\ref{noninteger3}), so we can confine to assume that there exists at least two such rays.

Assume first that there exists $j\in I$ such that $R_j$ is nef.
Since, by \cite[Theorem 1.1]{Wi-OnContr}, the general fiber of any contraction $\varphi_{R_j}$, $j\in I$, of fiber type has dimension $\geq 2(n-3)$, we derive that we can have at most two such contractions. Moreover, in this case $n=4$. It follows that both the contractions are as in case (5-3) of Proposition (\ref{noninteger3}); so $X=\pd\times\pd$ by \cite[Theorem A]{Sa-2Pbundle}, hence we are in case (1). 
\\
We claim that there cannot be any $i\in I$ such that $R_i$ is non nef. Indeed, if this is not the case, then the general fiber of the contraction $\varphi_{R_i}$ has dimension $\geq 2(n-3)+1$ by \cite[Theorem 1.1]{Wi-OnContr}; then, recalling that fibers of different extremal rays can meet only at points, we get a contradiction.

Assume now that all the $R_i$, $i\in I$, are birational. Then each $\varphi_{R_i}$ is as in case (5-4) of Proposition (\ref{noninteger3}) and the exceptional loci of the $R_i$'s are disjoint.
Furthermore, $(E_i,L_{E_i})=(\pt,\Ol_\pt(2))$ for any exceptional divisor $E_i$. So we are in case (2) of the statement.
\qed

\smallskip

In the following proposition we consider lower values of $\frac{1}{\sigma(X,L)}$.
%
\begin{proposition}
Let $(X,L)$ be a rays-positive manifold of dimension $n$. Assume that $0 \not= n-k < \frac{1}{\sigma(X,L)} < n-k+1$, $k \geq 4$ and $n \geq 2k-3$.
Then $(X,L)$ admits exactly one extremal ray satisfying $(\sigma(X,L) K_X+ L) \cdot R=0$ and it is described in Proposition (\ref{nonintegerk}), unless one of the following holds:
\begin{itemize}
	\item[(1)] $n=5$, $\sigma(X,L)=\frac{2}{3}$ and there exist extremal rays $R_1, \dots, R_i$, $i\in\{2, \dots, m\}$, such that $(2K_X+3L)\cdot R_i=0$; moreover, one of the following holds:   
\begin{itemize}
	\item[(1-1)] $m=2$ and each $\varphi_{R_i}$ is as in case (7-9) of Proposition (\ref{nonintegerk});
	\item[(1-2)] $m=2$, $\varphi_{R_1}$ is as in case (7-8) of Proposition (\ref{nonintegerk}) with $\varphi_{R_1}(X)=\pd$ and $\varphi_{R_2}$ is as in case (7-9) of Proposition (\ref{nonintegerk});
	\item[(1-3)] $\varphi_{R_1}$ is as in case (7-9) of Proposition (\ref{nonintegerk}) and each $\varphi_{R_i}$, with $i \geq 2$, is the blow-up of a smooth variety along a smooth curve such that $E_j \cap E_k=\emptyset$ if $j \not= k$ ($j,k \geq 2$);
\end{itemize}
	\item[(2)] $n=5$, $\sigma(X,L)=\frac{3}{4}$ and $X$ is the blow-up of a smooth variety at a finite set of points;
	\item[(3)] $n=2k-2$, $\sigma(X,L)=\frac{2}{n-1}$ and $X$ is the blow-up of a smooth variety at a finite set of points;
	\item[(4)] $n=2k-3$, $\sigma(X,L)=\frac{2}{n-2}$ and $X$ is the simultaneous contraction of extremal rays as in cases (7-11) and/or (7-12) of Proposition (\ref{nonintegerk}) with disjoint exceptional divisors.
\end{itemize}
\end{proposition}

\proof
Let $\{R_i:=\mathbb R_+ [\Gamma_i]\}_{i\in I}$ be the set of all the extremal rays in $\overline{\cone}(X)$ satisfying $(\sigma(X,L) K_X+ L) \cdot R_i=0$ and denote by $\tau:=\tau_L(R_i)=\frac{1}{\sigma(X,L)}$.
Then each $R_i$ satisfies the assumptions of Proposition (\ref{nonintegerk}), so we can confine to assume that there exists at least two such rays.

Assume first that there exists $j\in I$ such that $R_j$ is nef.
Since, by \cite[Theorem 1.1]{Wi-OnContr}, the general fiber of any contraction $\varphi_{R_j}$, $j\in I$, of fiber type has dimension $\geq 2(n-k)$, we derive that we can have at most two such contractions. Moreover, in this case $n=5$ and $k=4$. 
Now it is straightforward to get cases (1-1) and (1-2) of the statement by Proposition (\ref{nonintegerk}).
\\
Assume that there exists a non nef ray $R_i$ for some $i\in I$. In this case the general fiber of the contraction $\varphi_{R_i}$ has dimension $\geq 2(n-k)+1$ by \cite[Theorem 1.1]{Wi-OnContr}; again we have $n=5$ and $k=4$, so it is immediate to get case (1-3) of the statement by Proposition (\ref{nonintegerk}).

Assume now that all the $R_i$, $i\in I$, are birational. Then each $\varphi_{R_i}$ is as in one of cases (6-3), (7-10)--(7-12) of Proposition (\ref{nonintegerk}) and the exceptional loci of the $R_i$'s are disjoint.
Now it is straightforward to get cases (2)--(4) of the statement.
\qed

\subsection{Application to polarized manifolds}

In this section we apply the previous results to describe the nefvalue morphism of polarized manifolds.

\smallskip

Let $(X,L)$ be a {\em polarized manifold} of dimension $n\geq 3$.
We can define the {\em nefvalue} of $(X,L)$, which is the numerical invariant given by
$\tau := \tau (X,L) = \mbox{min} \{t \in \R: K_X+t L \mbox{ is nef}\}.$
Now, assume that $K_X$ is not nef, or, equivalently, that $X$ is not minimal in the sense of the {\em Minimal Model Program}, so that $\tau$ is a positive number.
By the {\em Kawamata's Rationality Theorem}, $\tau$ is a rational number.
Moreover, the divisor $K_X+\tau L$ defines a face $\Sigma := \{ C \in \overline\cone(X) : (K_X+\tau L) \cdot C = 0\}$ which is contained in the negative part (with respect to $K_X$) 
of the Kleiman--Mori cone 
and which is therefore generated by a finite number of extremal rays.
%
By the {\em Kawamata--Shokurov Base Point Free Theorem}, a high multiple of the divisor $K_X+\tau L$ is spanned by global sections, so it defines a morphism
$\Phi \colon X \to Y$ onto a normal projective variety with connected fibers. The map $\Phi$ is called the {\em nefvalue morphism} (relative to $(X,L)$).
Note that by construction $-K_X$ is $\Phi$-ample, therefore $\Phi$ is a Fano--Mori contraction 
and it contracts all curves in~$\Sigma$.

The nefvalue morphism is classically studied in the {\em Adjunction Theory} (see \cite{bookBS}) to describe polarized manifolds; however, by the discussion above, it is clear that it can be studied by looking at it as the contraction of an extremal face, which factors through the contraction of  extremal rays. 
Notice that, for any extremal ray $R$ orthogonal to $K_X+\tau(X,L)L$, we have $\tau_L(R)=\tau(X,L)=\frac{1}{\sigma(X,L)}$, where $\sigma(X,L)$ is the invariant defined in (\ref{defsigma}).
Moreover, since $n \geq 3$ and we are interested in non-integral values of $\tau_L(R)$, in view of \cite[Proposition 1.2]{BKLN-LPosRays}, we have $\tau_L(R)< n-1$.

The next proposition deals with polarized manifolds such that $n-2 < \tau(X,L) < n-1$ (for the proof in terms of Adjunction Theory, cf. \cite[Theorem 7.3.4]{bookBS}).
\begin{prop}
Let $(X,L)$ be a polarized manifold of dimension $n$. Assume that $K_X$ is not nef and that $n-2 < \tau(X,L) < n-1$.
Then the nefvalue morphism $\Phi$ is an elementary contraction and $(X,L)$ is described in \cite[Proposition 1.4]{BKLN-LPosRays}.
\end{prop}

\smallskip

In the following proposition we take care of the case $0 \not= n-3 < \tau(X,L) < n-2$ (for the proof in terms of Adjunction Theory, cf. \cite[Section 7]{bookBS}).
\begin{prop}
Let $(X,L)$ be a polarized manifold of dimension $n$. Assume that $K_X$ is not nef and that $0 \not= n-3 < \tau(X,L) < n-2$.
Then the nefvalue morphism $\Phi$ is an elementary contraction and $(X,L)$ is described in Proposition (\ref{noninteger3}), unless $\tau(X,L)=\frac{3}{2}$ and one of the following holds:
\begin{itemize}
	\item[(1)] $(X,L)=(\pd\times\pd, \Ol_{\pd\times\pd}(2,2))$ and $\Phi$ contracts $X$ to a point;
	\item[(2)] $\Phi\colon X \to Y$ is the blow-up of the smooth variety $Y$ at a finite number of points and $(E_i,L_{E_i})=(\pt,\Ol_\pt(2))$
	 for all exceptional divisors $E_i$.
\end{itemize}
\end{prop}

\smallskip

In the following proposition we consider lower value of $\tau(X,L)$ (for the proof in terms of Adjunction Theory, cf. \cite[Theorem 2.1]{BdT-nonint}).
\begin{prop}
Let $(X,L)$ be a polarized manifold of dimension $n$. Assume that $K_X$ is not nef and that $0 \not= n-k < \tau(X,L) < n-k+1$, $k \geq 4$ and $n \geq 2k-3$.
Then the nefvalue morphism $\Phi$ is an elementary contraction and $(X,L)$ is described in Proposition (\ref{nonintegerk}), unless one of the following holds:
\begin{itemize}
	\item[(1)] $n=5$, $\tau(X,L)=\frac{3}{2}$ and there exist extremal rays $R_1, \dots, R_i$, $i\in\{2, \dots, m\}$, such that $(K_X+\tau(X,L)L)\cdot R_i=0$; moreover, one of the following holds:   
\begin{itemize}
	\item[(1-1)] $m=2$ and each $\varphi_{R_i}$ is as in case (7-9) of Proposition (\ref{nonintegerk});
	\item[(1-2)] $m=2$, $\varphi_{R_1}$ is as in case (7-8) of Proposition (\ref{nonintegerk}) with $\varphi_{R_1}(X)=\pd$ and $\varphi_{R_2}$ is as in case (7-9) of Proposition (\ref{nonintegerk});
	\item[(1-3)] $\varphi_{R_1}$ is as in case (7-9) of Proposition (\ref{nonintegerk}) and each $\varphi_{R_i}$, with $i \geq 2$, is the blow-up of a smooth variety along a smooth curve such that $E_j \cap E_k=\emptyset$ if $j \not= k$ ($j,k \geq 2$);
\end{itemize}
	\item[(2)] $n=5$, $\tau(X,L)=\frac{4}{3}$ and $\Phi$ is the blow-up of a smooth variety at a finite set of points;
	\item[(3)] $n=2k-2$, $\tau(X,L)=\frac{n-1}{2}$ and $\Phi$ is the blow-up of a smooth variety at a finite set of points;
	\item[(4)] $n=2k-3$, $\tau(X,L)=\frac{n-2}{2}$ and $\Phi$ is the simultaneous contraction of extremal rays as in cases (7-11) and/or (7-12) of Proposition (\ref{nonintegerk}) with disjoint exceptional divisors.
\end{itemize}
\end{prop}

\bibliographystyle{plain}

\end{document}